\documentclass[a4paper,12pt,reqno]{amsart}
\usepackage[utf8]{inputenc}
\usepackage[T1]{fontenc}
\usepackage{lmodern}
\usepackage[english]{babel}
\usepackage{amsmath,a4wide}
\usepackage{xfrac}
\usepackage{esint}
\usepackage{mathrsfs,bm,amsthm,mathtools,yfonts,amssymb,color,braket,booktabs,graphicx,graphics,amsfonts}
\usepackage{latexsym,microtype,indentfirst,hyperref}
\usepackage{xcolor}
\usepackage{courier}
\usepackage[colorinlistoftodos]{todonotes}
\newcommand{\llbracket}{[\![}
\newcommand{\rrbracket}{]\!]}

\begingroup
\newtheorem{theorem}{Theorem}[section]

\endgroup

\theoremstyle{definition}
\newtheorem{definition}[theorem]{Definition}

\newtheorem{assumption}[theorem]{Assumption}

\newcommand\res{\mathop{\hbox{\vrule height 7pt width .3pt depth 0pt\vrule height .3pt width 5pt depth 0pt}}\nolimits}

\newcommand{\RR}{\mathbb{R}}

\newcommand{\HH}{\mathcal{H}}

\newcommand{\bB}{\mathbf{B}}
\usepackage{xcolor}

\def\a#1{\left\llbracket{#1}\right\rrbracket}

\title[Allard-type boundary regularity theory]{Allard-type regularity theory for area minimizing currents at boundaries with arbitrary multiplicity}

\author[I. Fleschler]{Ian Fleschler}
\address{Department of Mathematics, Fine Hall, Princeton University, Washington Road, Princeton, NJ 08540, USA}
\email{imf@princeton.edu}

\begin{document}
\begin{abstract}
This is an announcement of a series of upcoming works on boundary regularity for area minimizing currents, one of which is in collaboration with Reinaldo Resende. The setting we consider is that of an area minimizing current with a smooth boundary taken with arbitrary multiplicity. The main result is a generalization of Allard's boundary regularity theorem from which we derive important structural consequences.
\end{abstract}
\maketitle

\section{Introduction}
 This is an announcement of a series of upcoming works on boundary regularity for area minimizing currents. The setting we consider is that of an area minimizing current with a smooth boundary taken with an arbitrary multiplicity $Q$. We study the boundary points of minimum density (i.e., density $Q/2$). In the first work, by the named author, we show the uniqueness of the tangent cone for area minimizing currents at boundary points of density $Q/2$. These tangent cones are open books: they consist of a union of half-planes with multiplicities meeting at a common spine. In the second work, by the named author in collaboration with Reinaldo Resende, we show that the set of density $Q/2$ boundary singular points for an $m$-dimensional area minimizing current is $\HH^{m-3}$ rectifiable.  We also show that the boundary regular set (with no restriction on the density) is open and dense. In the third work, by the named author, we provide an example of a $3$ dimensional area minimizing current in $\RR^5$  with a boundary singularity of the type considered in our regularity theory. This example shows the optimality in the dimension of the above theorems.

We can also show that under the assumption of a convex barrier, which was previously considered by Allard, every boundary point is of density $Q/2$. Therefore, our theorems apply to this setting. Namely, if the boundary with multiplicity $Q$ lies on the boundary of a uniformly convex open set, then the boundary tangent cones are unique everywhere and the singular set is $\HH^{m-3}$ rectifiable.   

In particular, we settle the questions posed in Section 16.6 of De Lellis ICM Survey \cite{CamilloICM}.
\subsection{Setup}
The setting for the discussion of boundary regularity is the following:
\begin{assumption}[Assumption $C^{k,\alpha}$]\label{A:general}
Let integers $m, n\geq 2$, $Q\geq 1$, $k \in \mathbb{Z}_{\geq 0} \cup \left\{\infty, \omega \right\}$ and $\alpha \in [0,1]$. Consider $\Gamma$ a $C^{k,\alpha}$ oriented $(m-1)$-submanifold of $\RR^{m+n}$ without boundary, and assume $0\in\Gamma$ (when $k=\omega$, $\Gamma$ is a real-analytic submanifold of $\RR^{m+n}$). Let $T$ be an integral $m$-dimensional area minimizing current in $\bB_1$ with boundary $\partial T \res \bB_1=Q\a{\Gamma \cap \bB_1}$.
\end{assumption}
\newpage
\begin{definition} Let $T,\Gamma$ be as in Assumption \ref{A:general}. 
\begin{itemize}
    \item We say $0$ is a \textit{one-sided regular point} if there exists $\rho>0$ such that
\begin{equation*}
T \res \bB_{\rho}=\sum_{i=1}^N Q_i\a{\Sigma_i}
\end{equation*}
where $\Sigma_i$ are smooth minimal surfaces in $\bB_{\rho}$ with boundary $\Gamma$, disjoint outside of $\Gamma$, and meet transversally at $\Gamma$.
\item We say $0$ is a \textit{two-sided regular point} if for some $\rho>0$, there is a smooth minimal surface $M$ which contains $\Gamma$ and is being split into two pieces $M^+$ and $M^-$ by $\Gamma$ with
\begin{equation*}
T \res \bB_{\rho}=(Q+Q ^*)\a{M^{+}}+Q^*\a{M^{-}}
\end{equation*}
where $\partial \a{M^+}=-\partial \a{M^-}=\a{\Gamma}$.
\end{itemize}
\end{definition}
In these works we mainly study boundary points of minimum possible boundary density, which is instead guaranteed by the convex barrier condition:
\begin{assumption}[Convex Barrier] \label{convexbarrier} Let $\Omega \subset \RR^{m+n}$ be a domain such that $\partial \Omega$ is a $C^2$ uniformly convex submanifold of $\RR^{m+n}$. We say that $\sum_{i=1}^N Q_i \a{\Gamma_i}$ has a convex barrier if $Q_i$ are positive integers and $\Gamma_i \subset \partial \Omega$ are disjoint $C^2$ closed oriented submanifolds of $\partial \Omega$. In this setting we consider $T$ an area minimizing current with $\partial T=\sum_i Q_i \a{\Gamma_i}.$ 
\end{assumption}

\subsection{Aim of the works}
We are interested in further understanding the boundary regularity theory of $T$ when the boundary is taken with arbitrary multiplicity.

We begin by noticing that this problem can only have potentially interesting singular behaviour when the co-dimension is at least $2$. This is a consequence of  the following theorem of Hardt-Simon \cite{HS}. 
\begin{theorem}[Hardt-Simon \cite{HS}] Let $T$ be an area minimizing current which satisfies Assumption \ref{A:general} with $k=2, \alpha=0$, and codimension $n=1$. Then every boundary point is either a one-sided or a two-sided regular point.
\end{theorem}
The aim of these works is to study a generalization of the following results to the case $Q>1$:
\begin{theorem}[Allard Boundary Regularity \cite{AllB}]\label{t:allardb}
Let $T$ be an area minimizing current which satisfies Assumption \ref{A:general} with $k=2$, $\alpha=0$, $Q=1$ and assume $\Theta(T,0)=1/2$. Then $0$ is a one-sided regular point of $T$ (i.e. $T$ agrees with a classical smooth minimal surface with boundary at a neighborhood of $0$). Furthermore, if $\partial T=\sum_{i=1}^N Q_i \a{\Gamma_i}$ has a convex barrier (i.e. it satisfies Assumption \ref{convexbarrier}), then every point of $\Gamma$ is a one-sided regular point for $T$.
 \end{theorem}

 \begin{theorem}[Open and dense \cite{DDHM}]\label{t:fournames}
 Let $T$ be an area minimizing current which satisfies Assumption \ref{A:general} with $k=3$, $\alpha>0$ and $Q=1$. Then the boundary regular set (which consists of both one-sided and two-sided regular points) is open and dense. 
 
 Moreover, this result is optimal in the Hausdorff dimension: there exists a $2$d area minimizing current $T$ with a smooth boundary curve $\Gamma$ taken with multiplicity $Q=1$, such the boundary singular set (the complement of the regular set) has Hausdorff dimension $1$.
 \end{theorem}

 The program was successfully already carried in dimension $2$ by the following results:
 \begin{theorem}[Arbitrary multiplicity Allard-type Boundary Regularity  in dimension $2$ \cite{delellis2021allardtype}]
 Let $T$ be an area minimizing current which satisfies Assumption \ref{A:general} with $k=3$, $\alpha>0$, $m=2$ and assume $\Theta(T,0)=Q/2$. Then $0$ is a one-sided boundary regular point of $T$. Furthermore, if $\partial T=\sum_{i=1}^N Q_i \a{\Gamma_i}$ has a convex barrier (i.e. it satisfies Assumption \ref{convexbarrier}), then every point of $\bigcup_{i=1}^N \Gamma_i$ is a one-sided regular point for $T$.
 \end{theorem}

 \begin{theorem}[Open and dense \cite{reinaldo2024density}]
 Let $T$ be area minimizing which satisfies Assumption \ref{A:general} with $k=3, \alpha>0$ and $m=2$. Then the boundary regular set (which consists of both one-sided and two-sided regular points) is open and dense. As in Theorem \ref{t:fournames},
 this result is optimal in the Hausdorff dimension. \end{theorem}

\subsection{Main theorems}
In the works \cite{ian2024uniqueness,ianreinaldo2024regularity,ian2024example} we prove appropriate generalizations of Theorems \ref{t:allardb} and \ref{t:fournames} to higher multiplicity $Q>1$ and higher dimension $m \geq 3$:
\begin{theorem}[Rectifiability of the one-sided singular set \cite{ianreinaldo2024regularity}]\label{t:retifiability} Let $T$ be an area minimizing current which satisfies Assumption \ref{A:general} with $k=3$, $\alpha>0$ and $m \geq 3$. Further assume that $\Theta(T,p)=Q/2$ for every $p \in \bB_1 \cap \Gamma$. The boundary singular set (the set of points on $\bB_1 \cap \Gamma$ which are not one-sided regular points) is $\HH^{m-3}$ rectifiable (countable if $m=3$). Furthermore, if $\partial T=\sum_{i=1}^N Q_i \a{\Gamma_i}$ has a convex barrier (i.e. it satisfies Assumption \ref{convexbarrier}), the boundary singular set is $\HH^{m-3}$ rectifiable (countable if $m=3$). 
\end{theorem}

 \begin{theorem}[Open and dense \cite{ianreinaldo2024regularity}]
 Let $T$ be an area minimizing current which satisfies Assumption \ref{A:general} with $k=3$ and $\alpha>0$. Then the boundary regular set (which consists of both one-sided and two-sided regular points) is open and dense. As in Theorem \ref{t:fournames},
 this result is optimal in the Hausdorff dimension.
 \end{theorem}
 
 We also provide an example of a boundary singularity, which proves that Theorem 1.9 is dimensionally optimal. We remark that appropriate versions of the above theorems hold in ambient riemanian manifolds.

\subsection{Acknowledgments} The author acknowledges the support of the National Science Foundation through the grant FRG-1854147. Parts of this work took place during a two weeks visit of the author to Center of Nonlinear Analysis at Carnegie Mellon University and during a week visit of Reinaldo Resende to the Institute for Advanced Study.

The author is extremely grateful to his PhD advisor, Professor Camillo de Lellis, for his immense support, outstanding mentorship, and for proposing the above questions as the author’s PhD thesis topic.
The author also wishes to express gratitude for the fruitful collaboration to Reinaldo Resende, and to Paul Minter and Anna Skorobogatova for answering many relevant questions on these projects.
\bigskip

\section{Uniqueness of tangent cones} 
\subsection{History}
In 1993, Simon in \cite{simon1993cylindrical} studied the problem of uniqueness of tangent cones for a multiplicity $1$ class of stationary varifolds. In that paper, Simon introduced what are now known as Simon's estimates. The Simon's estimates use the remainder of the monotonicity formula in the presence of good enough density points - a condition typically referred to as a "no holes" condition - to prove that the $L^2$ height of the varifold is not concentrated at the spine of the cone. This allows a suitable reduction of the uniqueness of the tangent cone to understanding the rate of decay of energy for the linear problem.  
 
 The work of Simon was continued in a higher multiplicity setting in codimension $1$ by Wickramasekera in \cite{WickramasekeraSt} and applied again in both higher multiplicity and higher codimension in the context of area minimizing currents by De Lellis-Minter-Skorobogatova
\cite{de2023fineIII} and in parallel for the same setting by Krummel-Wickramasekera \cite{krummel2023analysisI,krummel2023analysisII}.

The recent paper \cite{DMSModp} also studies the problem in the context of area minimizing currents mod $p$ where the tangent cones considered are both open books (a union of planes meeting at a common spine) and planes. In the upcoming paper \cite{ian2024uniqueness}, the named author studies, for the first time, the problem in the context of an area minimizing current with boundary taken with a multiplicity $Q>1$ for an integral current of arbitrary dimension. We remark for experts that in this setting, the "no holes" condition is automatically satisfied at every point, which allows us to conclude uniqueness of the tangent cone at \textbf{every} minimum density boundary point. The main difference from \cite{DMSModp}
is that even though the class of cones contains open books like \cite{ian2024uniqueness} (they are those in the $m-1$ strata), the current does not have boundary mod $p$ and thus their work can be thought of instead as an interior regularity problem for a class of stationary varifolds. In addition, the class of cones they consider is more restricted, they must have a balancing condition: $\sum Q_i\nu_i=0$ where $\nu_i$ are the normals to the spine of the cone determined by the half-planes and $Q_i$ are their respective multiplicities,  which must also satisfy $Q_i \leq p/2$ and $\sum Q_i=p$. See also  \cite{Minter2024structure} for another reference of a uniqueness of tangent cones result with the class of tangent cones considered being a subclass of open books.

The situation and techniques are rather different in the context of $2$d area minimizing currents. In the paper \cite{White1983tangent}, White proves the uniqueness of the tangent cone for a $2$d area minimizing current using an epiperimetric inequality. In the paper \cite{hirsch2019uniqueness}, Hirsch-Marini adapt White's result to show the uniquenss of tangent cones at the boundary of $2$d area minimizing currents with boundary multiplicity $1$. Later, in the paper \cite{delellis2021uniqueness}, De Lellis-Nardulli-Steimbrüchel adapt the work of Hirsch-Marini for arbitrary boundary multiplicity.
\subsection{Results}
We start the program by obtaining a classification of the area minimizing cones. The boundary cones we consider are open books. 
\begin{definition}
A cone $C$ is an \textit{open book} if
 $C=\sum_{i=1}^N Q_i\a{H_i}$ such that  $H_i$ are distinct half-planes with  $\partial\a{H_i}=\a{V}$
where $V$ is an $m-1$ dimensional plane which is called the spine of the cone.  
\end{definition}

\begin{theorem}[Classification of the tangent cones \cite{ian2024uniqueness}]Let $C$ be an $m$-dimensional oriented area minimizing cone in $\RR^{m+n}$ with $\partial \a{C}=Q \a{\RR^{m-1}\times \left\{0\right\}}$.
Then $\Theta(C,0) \geq Q/2$. The equality case holds if and only if $C$ is an open book. Moreover there is a density jump: if $C$ is not an open book then $\Theta(C,0) \geq Q/2 +\varepsilon(m,n)$ for some $\varepsilon(m,n)>0$.
\end{theorem}

Under the Assumption \ref{A:general}
there is a monotonicity formula for the density at the boundary. As a consequence of the above classification, the boundary density must be at least $Q/2$ and equality only holds when every tangent cone is an open book. The density jump mentioned is only obtained at the end of the program as a consequence of the techniques for the main theorem rather than at the start.

Moreover, if $T$ is an area minimizing current with $\partial T=\sum_{i=1}^N Q_i \a{\Gamma_i}$ which satisfies the convex barrier Assumption \ref{convexbarrier}, then if $p \in \Gamma_i$, $\Theta(T,p)=Q_i/2$, and thus it belongs to the class of minimum density points, which we are studying.

The main theorem of the work is the following:
\begin{theorem}[Uniqueness of the tangent cone \cite{ian2024uniqueness}] 
Let $T$ be an area minimizing current which satisfies Assumption \ref{A:general} with $k=2,\alpha=0$. Further assume that 
$\Theta(T,0)=Q/2$. The tangent cone at $0$ is an open book and it is unique (i.e. it does not depend on the blow up sequence).
Furthermore, if $\partial T=\sum_{i=1}^N Q_i \a{\Gamma_i}$ has a convex barrier (i.e. it satisfies Assumption \ref{convexbarrier}), then $T$ has a unique tangent cone at every point of $\bigcup_{i=1}^N \Gamma_i$.
\end{theorem}
The previous statement in dimension $2$ was proved in \cite{delellis2021uniqueness}. The convex barrier case of the theorem follows from the general case, since every tangent cone arising from the convex barrier situation is an open book.

The most important consequence of the theorem above, which will enable us to run the regularity theory, is a decomposition theorem. This decomposition theorem allows us to reduce the study of boundary singularities to the study of flat boundary singularities. In dimension $2$, the decomposition was proved by \cite{delellis2021uniqueness,delellis2021allardtype}.
\begin{theorem}[Decomposition theorem \cite{ian2024uniqueness}]\label{t:decomposition} Let $T$ be an area minimizing current as in Assumption \ref{A:general} with $k=2$, $\alpha=0$. Further assume that $\Theta(T,0)=Q/2$. Let $C=\sum_{i=1}^N Q_i \a{H_i}$ be the unique tangent cone to $T$ at $0$ (where the representation of $C$ is such that the half-planes are distinct). Then there exists $\rho>0$ and area minimizing currents $T_1,T_2,...,T_N$ in $\bB_{\rho}$ such that
\begin{equation*}
T \res \bB_{\rho}=\sum_{i=1}^N T_i 
\end{equation*}
where the supports of $T_i$ only intersect at $\Gamma$, $\partial T_i \res \bB_{\rho}=Q_i \a{\Gamma}$, and the (unique) tangent cone at $0$ of $T_i$ is $Q_i \a{H_i}$.
\end{theorem}

\begin{section}{Rectifiabililty of the singular set}
\subsection{History} In the last decade Naber and Valtorta developed a flexible theory that makes it possible to prove the rectifiability of the singular set of geometric variational problems. Their theory uses the error term of the monotonicity formula to bound the Jones beta coefficients and reduces the original geometric variational problem to a problem in quantitative rectifiability. In the context of stationary varifolds see \cite{NV}. Their proof shows that each strata in Almgren's Stratification has finite measure, in particular the singular set of area minimizing currents in codimension $1$ needs to be $\HH^{m-7}$ rectifiable and of $\HH^{m-7}$ locally finite measure. The rectifiability of the singular set of area minimizing current in codimension $1$  was already known due to Simon \cite{simon1993rectifiability}, but his approach is much more rigid and does not show the finiteness of the measure. For relevant quantitative rectifiability results see the works of Azzam-Tolsa \cite{tolsa,azzamtolsa} and Edelen-Naber-Valtorta \cite{ENV}. For a more detailed discussion on the Naber-Valtorta methods along with their references and the question of $\sigma$-finiteness  of the measure see \cite{iancamillorectifiability}.

The question of rectifiability of the singular set in higher codimension was significantly more challenging. The Naber-Valtorta approach was implemented in that setting using the monotonicity formula for the frequency function, first for the 
for the linear problem in \cite{de2018rectifiability} and then later for the nonlinear setting by De Lellis-Skorobogatova in \cite{de2023fineI} and \cite{de2023fineII} to show that the points of frequency bigger than $1$ is $\HH^{m-2}$ rectifiable. This is a significant refinement of the works of De Lellis-Spadaro \cite{DS1,DS2,DS3,DS4,DS5} that build on Almgren's celebrated interior regularity paper \cite{Alm}. Together with \cite{de2023fineIII} De Lellis-Skorobogatova-Minter conclude the $\HH^{m-2}$ rectifiability of the singular set.

In the upcoming work \cite{ianreinaldo2024regularity}, which is in collaboration with Resende, we implement the Naber-Valtorta techniques in the boundary setting.

It should be noted that Krummel-Wickramasekera have carried out Simon's program to prove the rectifiability of the interior singular set first for the linear problem in \cite{krummel2017fine} and for the nonlinear problem in the works \cite{krummel2023analysisI,krummel2023analysisII} and an upcoming paper using the techniques mentioned in the previous section. Their work still requires the construction of the center manifold (cf. Chapter 4 \cite{Alm}, \cite{DS4}) and moreover an important starting point for \cite{krummel2017fine} is the uniqueness of tangent functions to $2$-dimensional Dir-minimizing Q-valued functions \cite{DS1}.

\subsection{Results} The results are the following:
\begin{theorem}[Rectifiability of the one-sided singular set \cite{ianreinaldo2024regularity}] Let $T$ be an area minimizing current which satisfies Assumption \ref{A:general} with $k=3$, $\alpha>0$ and $m \geq 3$. Further assume that $\Theta(T,p)=Q/2$ for every $p \in \bB_1 \cap \Gamma$. The boundary singular set (the set of points on $\bB_1 \cap \Gamma$ which are not one-sided regular points) is $\HH^{m-3}$ rectifiable (countable if $m=3$). Furthermore, if $\partial T=\sum_{i=1}^N Q_i \a{\Gamma_i}$ has a convex barrier (i.e. it satisfies Assumption \ref{convexbarrier}), the boundary singular set is $\HH^{m-3}$ rectifiable (countable if $m=3$). 
\end{theorem}
We remark that in this context even the dimensional bound was not known prior to our work in contrast to the interior case. The machinery of De Lellis-Spadaro \cite{DS1,DS2,DS3,DS4,DS5} was implemented in the boundary settings in \cite{DDHM,delellis2021allardtype,reinaldo2024density}. The dimensional bound is not a straightforward consequence of the above. Our techniques rely strongly on the finer analysis of the frequency function developed by De Lellis-Skorobogatova \cite{de2023fineI,de2023fineII}.
 \begin{theorem}[Open and dense \cite{ianreinaldo2024regularity}]
 Let $T$ be an area minimizing current which satisfies Assumption \ref{A:general} with $k=3$ and $\alpha>0$. Then the boundary regular set (which consists of both one-sided and two-sided regular points) is open and dense. As in Theorem \ref{t:fournames},
 this result is optimal in the Hausdorff dimension.
 \end{theorem}
Prior to the above theorem not even the existence of a single boundary regular point was known in this setting.
\end{section}
\begin{section}{Example}
\begin{definition} Let $T, \Gamma$ be as in Assumption \ref{A:general}. We say that $0$ is a \textit{removable one-sided singularity} if $\Theta(T,0)=Q/2$ and there exists $\rho>0$ such that
\begin{equation*}
T \res \bB_{\rho}=\sum Q_i \a{\Sigma_i}
\end{equation*}
where $\Sigma_i$ are smooth minimal surfaces and $\partial \Sigma_i=\Gamma \cap \bB_{\rho}$ (but they are not required to be disjoint outside of the boundary or meet transversal at the boundary). 
A point is an \textit{essential one-sided singularity} if $\Theta(T,0)=Q/2$, $0$ is a singular point but it is not a removable singularity.
\end{definition}
\begin{theorem}[Example \cite{ian2024example}] There exists a $3$d area minimizing current $T$ in $\bB_1$  with $T, \Gamma$ which satisfy Assumption \ref{A:general} with $k=\omega$, codimension $n=2$ (i.e. $\Gamma$ is a real-analytic submanifold of $\RR^5$) and $Q=2$ such that: $\Theta(T,p)=1$ for every $p \in \Gamma \cap \bB_1$ and $0$ is an essential one-sided boundary singularity.
\end{theorem}
In particular, the consequence of this is that the dimensional bound of the one-sided singular set of Theorem \ref{t:retifiability} is sharp in any dimension, even when the boundary is required to be real-analytic.
\end{section}
\begin{section}{Linear problem}
The results outlined in the previous sections also have analogous results in the setting of the boundary linear problem. 

\begin{assumption}
Given a Dir-minimizing function $u \in W^{1,2}(\Omega,A_{Q}(\RR^n))$ for an open set $\Omega \subset \RR^m$. The boundary linear problem consists of functions $u$ such that $u \res \partial \Omega \cap \bB_1=Q\a{0}$.
\end{assumption}
We also have theorems analogous to those presented above under suitable assumptions on the regularity of the domain, namely:
\begin{enumerate}
\item The frequency function is at least $1$ at the boundary. The only $1$- homogeneous blowups are linear and the blowups are unique at every frequency $1$ point. Moreover there is a frequency jump (i.e. if the frequency is not $1$ then it must be above $1+\varepsilon(m,n)$ for some $\varepsilon(m,n)>0$).
\item The boundary singular set for the linear problem is $\HH^{m-3}$ rectifiable.
\item There is a $2$-valued $I>1$ homogeneous Dir minimizer in $W^{1,2}(\RR^2 \times [0,\infty),A_{2}(\RR^2))$  with zero boundary value.  In $2$d this is not possible since the only homogeneous Dir minimizers are linear (and thus $I=1$ is the only option).
\end{enumerate}
\end{section}
 \bibliographystyle{alpha}
        \bibliography{Cortona}

\end{document}